\documentclass[authoryear,12pt]{elsarticle}
\usepackage{caption}
\usepackage{color}
\usepackage{lscape}
\usepackage{afterpage}
\usepackage{pstricks}
\usepackage{pst-plot}
\usepackage{longtable}
\usepackage{dcolumn}
\usepackage{pst-node}
\usepackage{amsmath}
\usepackage{amssymb}
\usepackage{amsthm}
\usepackage{multirow}
\usepackage{colortab}
\usepackage{color}
\usepackage{array}
\usepackage{slashbox}
\usepackage{colortbl}
\usepackage{subfigure}
\usepackage{textcomp}
\usepackage{pstricks, pst-node}
\usepackage{pst-all}
\usepackage{algorithm,algorithmic}
\usepackage{url}
\usepackage{hyperref}
\psset{arrows=->, labelsep=3pt, mnode=circle}

\usepackage{eurosym}

\usepackage{tikz}

\newfont{\rams}{msbm10 scaled\magstep1}
\newcommand{\rea}{\mbox{\rams \symbol{'122}}}

\newenvironment{resumeT}{\begin{list}{}{\setlength{\rightmargin}{\leftmargin}}\item[]
{\centering {\bf \it~~~}
\par}\item[]\ignorespaces}{\unskip\end{list}}

\newtheorem{note}{Note}[section]

\begin{document}
\title{Scoring from Pairwise Winning Indices}

\author[Eco]{\rm Sally Giuseppe Arcidiacono}
\ead{s.arcidiacono@unict.it}
\author[Eco]{\rm Salvatore Corrente}
\ead{salvatore.corrente@unict.it}
\author[Eco,por]{\rm Salvatore Greco}
\ead{salgreco@unict.it}

\address[Eco]{Department of Economics and Business, University of Catania, Corso Italia, 55, 95129  Catania, Italy}
\address[por]{University of Portsmouth, Portsmouth Business School, Centre of Operations Research and Logistics (CORL), Richmond Building, Portland Street, Portsmouth PO1 3DE, United Kingdom}

\date{}
\maketitle

\vspace{-1cm}

\begin{resumeT}
{\large {\bf Abstract:}} The pairwise winning indices, computed in the Stochastic Multicriteria Acceptability Analysis, give the probability with which an alternative is preferred to another taking into account all the instances of the assumed preference model compatible with the information provided by the Decision Maker in terms of pairwise preference comparisons of reference alternatives. In this paper we present a new scoring method assigning a value to each alternative summarizing the results of the pairwise winning indices. Several procedures assigning a score to each alternative on the basis of the pairwise winning indices have been provided in literature. However, while all of them compute this score just to rank the alternatives under consideration, our method, expressing the score in terms of an additive value function, permits to disaggregate the overall evaluation of each alternative in the sum of contributions of considered criteria. This will permit not only to rank the alternatives but also to explain the reasons for which an alternative obtains its evaluation and, consequently, fills a certain ranking position. To prove the efficiency of the method in representing the preferences of the Decision Maker, we performed an extensive set of simulations varying the number of alternatives and criteria. The results of the simulations, analyzed from a statistical point of view, show the goodness of our procedure. The applicability of the method to decision making problems is explained by means of a case study related to the evaluation of financial funds. 


\vspace{0,3cm}
\end{resumeT}

\pagenumbering{arabic}

\section{Introduction}
Given a set of alternatives $A=\{a,b,\ldots\}$ evaluated on a coherent family of criteria \citep{Roy1996} $G=\{g_1,\ldots,g_m\}$, choice, ranking and sorting are the typical decision problems handled through Multiple Criteria Decision Making (MCDM) methods \citep{GrecoEhrgottFigueira2016,KeeneyRaiffa1976}. In this paper we are interested in ranking decision problems in which alternatives have to be ordered from the best to the worst taking into account the preferences of the Decision Maker (DM). Several methods aiming to produce such a ranking have been proposed in literature and they mainly differ (i) with respect to the form the DM's preferences are articulated and, (ii) in the procedures used to get the final ranking. With respect to the first point we distinguish between direct and indirect preference information. In case of direct preference information, the DM is asked to provide directly values of the parameters involved in the decision model used to produce the ranking, while, in case of indirect preference information, the DM is asked to provide some information in terms of preference comparison between alternatives or comparison of criteria with respect to their importance; this information is then used to infer parameters of the assumed preference model so that the application of the model with the inferred parameters restores the preferences expressed by the DM. The indirect preference information is the most used in practice because of the lower cognitive effort asked to the DM in providing it. The approach applying this way of asking preferences to the DM is known as preference disaggregation \citep{JacquetLagrezeSiskos1982,JacquetLagrezeSiskos2001}. 

As explained above, in preference disaggregation one aims to infer an instance of the assumed preference model compatible with the preferences provided by the DM. However, many of these instances can exist and, while all of them give the same recommendation on the comparison of alternatives on which the DM expressed her preferences, they can provide different recommendations comparing other alternatives. For such a reason, giving a recommendation using only one of these instances can be considered arbitrary to some extent. Robust Ordinal Regression (ROR) \citep{CorrenteEtAl2013,GrecoMousseauSlowinski2008} and Stochastic Multicriteria Acceptability Analysis (SMAA) \citep{LahdelmaHokkanenSalminen1998,PelissariEtAl2020} are two families of MCDM methods aiming to provide recommendations on the problem at hand taking into account all the instances of the preference model compatible with the preferences given by the DM. On the one hand, ROR builds necessary and possible preference relations for which $a$ is necessarily preferred to $b$ iff $a$ is at least as good as $b$ for all compatible models, while $a$ is possibly preferred to $b$ iff $a$ is at least as good as $b$ for at least one compatible model. On the other hand, SMAA gives recommendations in statistical terms producing mainly two indices: (i) the Rank Acceptability Index, $b^{k}(a)$, being the probability with which an alternative $a$ fills the position $k$ in the final ranking, (ii) the Pairwise Winning Index (PWI), $p(a,b)$, being the probability with which $a$ is at least as good as $b$. In this paper we are interested in the SMAA methodology and, in particular, in the aggregation of the PWIs to get a complete ranking of the alternatives under consideration. Several methods have been proposed in literature to deal with this problem. In the following, without any ambition to be exhaustive, we shall review some of them. The first methods appear in social choice theory \citep{ArrowSenSuzumura2010}. For example,  \cite{Dodgson1876} proposes a ranking method such that best ranked alternatives are the closest to be Condorcet winner \citep{Condorcet1785}. This is the alternative $a\in A$ (if it exists) such that, with respect to any other alternative $b\in A$ there is a majority of compatible value functions for which $a$ is preferred to $b$ (that is, the alternative $a\in A$ for which $p(a,b)\geqslant k$ for all $b\in A$, with $k\in]0.5,1]$ being a proper majority threshold) and, therefore, the alternative having a majority of preferences. Another very well-known ranking algorithm relating social choice with probability values $p(a,b)$ is the Simpson procedure \citep{Simpson1969} that ranks alternatives in $A$ assigning to each $a\in A$ a score being its minimal $p(a,b)$ over all $b$ in $A$. In the perspective of applying some social choice ranking procedure to PWIs, \cite{KadzinskiMichalski2016} present nine methods summarizing the PWIs by sum, min and max operators. In a more computational oriented approach, \cite{Vetschera2017} presents different ranking methods based on the solution of specific Mixed Integer Linear Programming (MILP) problems. 

Even if, as observed above, several methods have been proposed in the past to build a complete ranking of the alternatives summarizing the information provided by the PWIs, none of them provides any easily understandable explanation of this ranking. We then aim to fill this gap proposing a method that, on the basis of the PWIs supplied by the SMAA methodology gives a ranking of the considered alternatives and explains it through an additive value function. In this way, differently from the other methods proposed up to now in literature to obtain a final ranking from the PWIs of SMAA, the DM gets also a scoring for the considered alternatives and, moreover, using the value function provided by our procedure, she can investigate on the reasons for which an alternative gets a certain ranking position. Indeed, since the ranking is provided by means of an additive value function, the DM can look at the contribution given by each criterion to the global value permitting to identify criteria being weak and strong points for the considered alternatives. From a methodological point of view, as it will be clear later, the computation of the additive value function involves only to solve an LP problem. \\
Let us also observe that, in line with what has been proposed for ROR in \cite{KadzinskiGrecoSlowinski2012b}, the value function provided by our method can be seen as a value function representative of the many compatible value functions representing the preferences expressed by the DM. Therefore, our approach can be interpreted also as a procedure to construct a representative value function for SMAA methodology.

From a preference learning perspective \citep{FurkranzHullermeier2010a}, to prove the efficiency of our method not only in explaining the preferences of the DM but also in predicting her comprehensive preferences starting from some available information, we performed an extensive set of simulations considering different numbers of alternatives and criteria. With this aim, for each (\# Alternatives, \# Criteria) configuration, we assumed the existence of an artificial DM whose preferences, obtained through a value function that is unknown to the method we use to produce the representative value function, have to be discovered on the basis of some pairwise comparisons it provided. To evaluate the performances of our method with respect to this objective, we compute the Kendall-Tau correlation coefficient \citep{Kendall1938} between the ranking of the alternatives produced by the artificial DM and the one given by our procedure. We compared our method to other fourteen methods known in literature aiming to summarize the information contained in the PWIs. The results prove that there is not a method being the best for each (\# Alternatives, \# Criteria) configuration and the difference in terms of Kendall-Tau between the values obtained by our proposal and those obtained by the best method in each configuration is not statistically significant with respect to a Mann-Whitney U test with 5\% significance. This proves that our method is good not only for its capacity to \textit{explain} the preferences of the DM but also to learn preferences with reliable results. 

Finally, to show how the method can support decision making in a real world problem, we applied it to a financial context in which seven funds have to be overall evaluated taking into account their performance on five criteria. 

The paper is structured as follows. In the next section, we present the new method as well as some extensions; an extensive comparison between our proposal and other methods presented in literature to deal with the same problem is performed in Section \ref{ComparisonSection}, while the method is applied to a real world financial problem in Section \ref{CaseStudySection}; finally, some conclusions and further directions of research are given in Section \ref{ConclusionsSection}.

\section{Giving a score to the alternatives on the basis of the Pairwise Winning Indices}\label{ScoringSection}
\subsection{Methodological Background}
In the following, we shall briefly recall the main concepts of the preference disaggregation and SMAA. Without loss of generality, we shall assume that all criteria are of the gain type (the greater the evaluation of an alternative $a\in A$ on criterion $g_j\in G$, the more $a$ is preferred on $g_j$) and $g_j(a)$ will denote the evaluation of $a$ on $g_j$. Let us suppose that the model assumed to represent the preferences of the DM is an additive value function $U:A\rightarrow\left[0,1\right]$ of the following type

\begin{equation}\label{AdditiveVF}
U(a)=U\left(g_1(a),\ldots,g_m(a)\right)=\displaystyle\sum_{j=1}^{m}u_j\left(g_j(a)\right)
\end{equation}

\noindent where $u_j:X_j\rightarrow\left[0,1\right]$, for all $g_j\in G$, is the marginal value function related to $g_j$ and $X_{j}=\left\{x_j^{0},\ldots,x_j^{n_j}\right\}$ is the set of evaluations taken by alternatives in $A$ on $g_j$ such that $x_j^k<x_j^{k+1}$ for all $k=0,\ldots,n_{j}-1$. Moreover, $u_j$ is non-decreasing in $X_j$ for all $g_j\in G.$

\noindent To build an additive value function, the marginal value functions should be defined and, under the preference disaggregation approach, this is done by taking into account some preferences provided by the DM in terms of comparison between alternatives of the following type: 
\begin{itemize}
\item $a$ is preferred to $b$ (denoted by $a\succ_{DM}b$) translated into the constraint $U(a)>U(b)$,
\item $a$ is at least as good as $b$ ($a\succsim_{DM}b$) translated into the constraint $U(a)\geqslant U(b)$,
\item $a$ is indifferent to $b$ ($a\sim_{DM}b$) translated into the constraint $U(a)=U(b)$\footnote{The indifference between alternatives $a$ and $b$ can also be translated into a different constraint by using an auxiliary threshold as described in \cite{BrankeEtAl2017}. However, the way the indifference relation is translated does not affect the following description and, for this reason, we shall assume that $a\sim_{DM}b$ iff $U(a)=U(b)$.}.
\end{itemize}

A \textit{compatible value function} is then an additive value function of type (\ref{AdditiveVF}) compatible with the preferences provided by the DM and, therefore, satisfying the following set of constraints:

$$
\left.
\begin{array}{l}
U(a)>U(b), \;\;\mbox{iff}\; a\succ_{DM} b, \\[0,2cm]
U(a)\geqslant U(b), \;\;\mbox{iff}\; a\succsim_{DM} b, \\[0,2cm] 
U(a)=U(b), \;\;\mbox{if}\; a\sim_{DM} b, \\[0,2cm]
u_j\left(x_j^{k}\right)\leqslant u_j\left(x_j^{k+1}\right), \;\forall g_j\in G\;\mbox{and}\; k=0,\ldots,n_j-1,\\[0,2cm]
u_j\left(x_j^{0}\right)=0, \;\forall g_j\in G, \\[0,2cm]
\displaystyle\sum_{j=1}^{m}u_j\left(x_j^{n_j}\right)=1.
\end{array}
\right\}E^{DM}
$$

To check for the existence of at least one compatible value function one needs to solve the following LP problem 

$$
\begin{array}{l}
\;\;\varepsilon^*=\max\varepsilon\;\mbox{subject to} \\[0,2cm]
\left.
\begin{array}{l}
U(a)\geqslant U(b)+\varepsilon, \;\;\mbox{iff}\; a\succ_{DM} b, \\[0,2cm]
U(a)\geqslant U(b), \;\;\mbox{iff}\; a\succsim_{DM} b, \\[0,2cm] 
U(a)=U(b), \;\;\mbox{if}\; a\sim_{DM} b, \\[0,2cm]
u_j\left(x_j^{k}\right)\leqslant u_j\left(x_j^{k+1}\right), \;\forall g_j\in G\;\mbox{and}\; k=0,\ldots,n_j-1,\\[0,2cm]
u_j\left(x_j^{0}\right)=0, \;\forall g_j\in G, \\[0,2cm]
\displaystyle\sum_{j=1}^{m}u_j\left(x_j^{n_j}\right)=1
\end{array}
\right\}E^{DM'}
\end{array}
$$

\noindent where $\varepsilon$ is an auxiliary variable used to convert the strict inequalities ($U(a)>U(b)$) in weak ones ($U(a)\geqslant U(b)+\varepsilon$). If $E^{DM'}$ is feasible and $\varepsilon^*>0$, then, there exists at least one compatible value function. In the opposite case ($E^{DM'}$ is infeasible or $\varepsilon^*\leqslant 0$), then, it does not exist any compatible value function and the causes can be investigated by means of the methods presented in \cite{MousseauEtAl2003}. Let us assume that at least one compatible value function exists. In this case, as already observed in the introduction, in general, infinitely many compatible value functions exist and, therefore, the SMAA methodology aims to give a recommendation on the problem at hand by taking into account a sampling of them. Since the constraints in $E^{DM'}$ define a convex set of parameters, some compatible value functions can be sampled by using, for example, the Hit-And-Run algorithm \citep{Smith1984,TervonenEtAl2013,VanValkenhoefTervonenPostmus2014}. Let us denote by ${\cal U}$ the set of sampled compatible value functions. Each value function in ${\cal U}$ will give a certain recommendation on each pair of alternatives $(a,b)\in A\times A$. For this reason, for each $(a,b)\in A\times A$, two different subsets of ${\cal U}$ can be defined:
\begin{equation}\label{SubsetsValueFunctions}
{\cal U}_{a\succ b}=\{U\in {\cal U}: U(a)>U(b)\}; \;\;\; {\cal U}_{a\sim b}=\{U\in {\cal U}: U(a)=U(b)\}.
\end{equation}
For each $(a,b)\in A\times A$, the PWI of the pair $(a,b)$, denoted by $p(a,b)$, can then be computed in the following way: 
\begin{equation}\label{PWI}
p(a,b)=\frac{|{\cal U}_{a\succ b}|+\frac{1}{2}|{\cal U}_{a\sim b}|}{|{\cal U}|}.
\end{equation}
Let us observe that other definitions of the PWI have been provided in literature. For example, in \cite{KadzinskiMichalski2016}, $p(a,b)=\frac{|{\cal U}_{a\succsim b}|}{|{\cal U}|}$ where ${\cal U}_{a\succsim b}=\{U\in {\cal U}: U(a)\geqslant U(b)\}$ and, therefore, $\left|{\cal U}_{a\succsim b}\right|=\left|{\cal U}_{a\succ b}\right|+\left|{\cal U}_{a\sim b}\right|$, while in \cite{Vetschera2017} $p(a,b)=\frac{|{\cal U}_{a\succ b}|}{|{\cal U}|}$. In the following, without loss of generality, we shall consider the PWI defined in eq. (\ref{PWI}).

\subsection{Our proposal}\label{ScoringMethodSection}
The idea under the construction of a compatible value function able to represent in the best way the PWIs is the following: 

\begin{center}
``\textit{Given $a,b\in A$, if $p(a,b)\geqslant 0.5$, that is, for at least 50\% of the sampled compatible value functions $a$ is at least as good as $b$, then, the greater $p(a,b)$, the larger should be the difference between the utilities of $a$ and $b$.}" 
\end{center}

\noindent Given in different terms, we aim to build a compatible value function $U$ such that the difference $U(a)-U(b)$ increases with $p(a,b)$. Formally, this requirement is translated into the following constraint: 

\begin{equation}\label{ScoringConstraint}
\mbox{if}\;\;p(a,b)\geqslant 0.5, \;\;\mbox{then}\;\;U(a)-U(b)\geqslant\eta\left(p(a,b)-0.5\right)
\end{equation}

\noindent where $\eta\in\rea^{+}$. 

Let us observe that if $\eta>0$, the constraint (\ref{ScoringConstraint}) perfectly represents the eventual preferences provided by the DM on some pairs of alternatives. Indeed, 
\begin{itemize}
\item if $a\succsim_{DM}b$, then $p(a,b)\geqslant 0.5$ and, therefore, $U(a)-U(b)\geqslant\eta\cdot\left(p(a,b)-0.5\right)\geqslant 0$, from which it follows that $U(a)\geqslant U(b)$,
\item if $a\succ_{DM}b$, then $p(a,b)=1$ and, therefore, $U(a)-U(b)\geqslant 0.5\cdot\eta>0$ from which it follows that $U(a)>U(b)$, 
\item if $a\sim_{DM}b$, then $p(a,b)=0.5=p(b,a)$ and, therefore, on the one hand, $U(a)-U(b)\geqslant 0$ and, on the other hand, $U(b)-U(a)\geqslant 0$, so that we get $U(a)=U(b)$.
\end{itemize}
To check for an additive value function having all the characteristics mentioned above (we shall call it \textit{compatible scoring function}), the following LP problem, denoted by $LP_{0}$, should be solved: 
$$
\begin{array}{l}
\;\;\eta^*=\max \eta, \;\mbox{subject to}, \\[0,2cm]
\left.
\begin{array}{l}
U(a)-U(b)\geqslant \eta\cdot\left(p(a,b)-0.5\right), \;\forall (a,b)\in A\times A: p(a,b)\geqslant 0.5,\\[0,2cm]
u_j\left(x_j^{k}\right)\leqslant u_j\left(x_j^{k+1}\right), \;\forall g_j\in G\;\mbox{and}\;k=0,\ldots,n_j-1,\\[0,2cm]
u_j\left(x_j^{0}\right)=0, \;\forall g_j\in G, \\[0,2cm]
\displaystyle\sum_{j=1}^{m}u_j\left(x_j^{n_j}\right)=1.
\end{array}
\right\}E^{SF}
\end{array}
$$

\noindent If $E^{SF}$ is feasible and $\eta^{*}>0$, then there is at least one compatible scoring function. In the opposite case ($E^{SF}$ is infeasible or $\eta^{*}\leqslant 0$), then there is not any compatible value function satisfying all constraints in $E^{SF}$ with a positive value of $\eta$. The causes of this infeasibility could be detected by using one of the approaches proposed in \cite{MousseauEtAl2003}. In the first case ($E^{SF}$ is feasible and $\eta^{*}>0$), the compatible scoring function obtained by solving $LP_{0}$ can then be used to assign a value to each alternative in $A$ producing, therefore, a complete ranking of the alternatives at hand.  Observe that the value assigned to the alternatives can be considered also as a scoring obtained taking into account the whole set of compatible value functions that, in fact, have been considered to compute the PWIs. In this perspective the obtained compatible scoring function can be seen as a representative value function \citep{KadzinskiGrecoSlowinski2012a} built on the basis of the principle ``one for all, all for one". Indeed, on the one hand, all compatible value functions concur to define the value function $U$ corresponding to $\eta^\ast$, while, on the other hand, $U$ represents the whole set of the compatible value functions. 

Let us observe that, as will be described more in detail in Section \ref{ComparisonSection}, even if $E^{SF}$ is feasible but $\eta^*\leqslant 0$, the compatible scoring function obtained solving $LP_{0}$ can be used to assign a value to each alternative and, then, ranking completely all alternatives under consideration. In this case, there will be at least one pair of alternatives $(a,b)\in A\times A$ for which $p(a,b) \geqslant 0.5$ and nevertheless $U(a)<U(b)$. However, maximizing $\eta$, that is, in this case, minimizing its absolute value since it is negative, the optimization problem gives a value function $U$ which minimizes the deviation from the preferences represented by the PWIs. 

\subsection{Some extensions of our scoring method}\label{ExtensionsSection}
Let us assume that $E^{SF}$ is feasible and $\eta^{*}>0$. This means that at least one compatible scoring function exists. However, many of them could exist. Let us denote by ${\cal U}^{SF}$ the set of additive value functions satisfying all constraints in $E^{SF}$ with $\eta=\eta^{*}$ that is, all value functions maximally discriminating between the alternatives $(a,b)\in A\times A$ on the basis of the corresponding PWIs. From now on, we shall call such functions \textit{maximally discriminant compatible scoring functions}. As already observed in the introduction, the previous scoring procedures based on PWIs compute a value for each alternative just to obtain a complete ranking of them. However, this number has not a particular meaning and, in some cases, can be completely useless for the DM. In our approach, the complete ranking of the alternatives under examination is produced on the basis of the construction of a compatible scoring function that, on the one hand, assigns a score to each alternative and, on the other hand, \textit{explains} the reasons for which an alternative has received a particular ranking position. This explanation is provided defining the contribution of each criterion $g_j\in G$ to the overall evaluation through the values taken by the corresponding marginal value function $u_j$.  For this reason, from the explainability point of view (see e.g. \citealt{ArrietaEtAl2020}), it is important checking if among the maximally discriminant compatible scoring functions in ${\cal U}^{SF}$, there exists at least one having one of the following characteristics: 
\begin{description}
\item[(T1)] all criteria $g_j\in G$ give a contribution to the compatible scoring function $U$,
\item[(T2)] for each $g_j\in G$ the marginal value function $u_j$ is strictly monotone in $X_j$.
\end{description}

Considering \textbf{(T1)}, one is looking for a maximally discriminant compatible scoring function $U\in{\cal U}^{SF}$ such that $u_j\left(x_{j}^{n_j}\right)>0$ for each $g_j\in G$. Indeed, in an additive value function, $u_j\left(x_{j}^{n_j}\right)$ can be considered, in some way, as the ``importance" of criterion $g_j$ being the marginal value assigned to the greatest performance on $g_j$ (that is $x_{j}^{n_j}$) by the marginal value function $u_j$ \citep{JacquetLagrezeSiskos1982}.   To this aim, the following LP problem, denoted by $LP_1$, has to be solved: 
$$
\begin{array}{l}
\;\;h^*=\max\; h, \;\mbox{subject to}, \\[0,2cm]
\left.
\begin{array}{l}
\eta=\eta^{*}, \\[0,2cm]
u_j\left(x_j^{n_j}\right)\geqslant h, \;\;\forall g_j\in G, \\[0,2cm]
E^{SF}.
\end{array}
\right\} E_{AllContr}^{SF}
\end{array}
$$ 
If $E_{AllContr}^{SF}$ is feasible and $h^{*}>0$, then, there is at least one maximally discriminant compatible scoring function $U\in{\cal U}^{SF}$ in which all criteria contribute to $U$. In the opposite case, in all maximally discriminant compatible scoring functions, $u_j\left(x_j^{n_j}\right)=0$ and, therefore, $g_j$ is not contributing to the global value assigned to the alternatives by $U$. Let us denote by ${\cal U}^{SF}_{AllContr}$ the subset of ${\cal U}^{SF}$ composed of all maximally discriminant compatible scoring functions in which all criteria contribute to the global value of each alternative.

Considering \textbf{(T2)}, one is looking for a maximally discriminant compatible scoring function $U\in{\cal U}^{SF}$ where, for each $g_j\in G$, $u_j\left(x_j^{k}\right)<u_j\left(x_j^{k+1}\right)$ for all $k=0,\ldots,n_j-1$. To check for the existence of a maximally discriminant compatible scoring function having the mentioned characteristic, one has to solve the following LP problem denoted by $LP_2$:
$$
\begin{array}{l}
\;\;\sigma^*=\max\; \sigma, \;\mbox{subject to}, \\[0,2cm]
\left.
\begin{array}{l}
\eta=\eta^{*}, \\[0,2cm]
u_j\left(x_j^{k}\right)+\sigma\leqslant u_j\left(x_j^{k+1}\right), \;\;\forall g_j\in G \;\mbox{and}\;k=0,\ldots,n_j-1, \\[0,2cm]
E^{SF}.
\end{array}
\right\} E_{AllInc}^{SF}
\end{array}
$$
If $E_{AllInc}^{SF}$ is feasible and $\sigma^{*}>0$ then there exists at least one maximally discriminant compatible scoring function in ${\cal U}^{SF}$ for which all marginal value functions are increasing, while, in the opposite case ($E_{AllInc}^{SF}$ is infeasible or $\sigma^{*}\leqslant 0$), then this is not the case. Let us denote by ${\cal U}^{SF}_{AllInc}$ the subset of ${\cal U}^{SF}$ composed of the maximally discriminant compatible scoring functions with $\eta=\eta^*$ such that all marginal value functions are increasing. 

\begin{note}\label{InclusionNote}
Let us observe that ${\cal U}^{SF}_{AllInc}\subseteq{\cal U}^{SF}_{AllContr}$ and, therefore, if ${\cal U}^{SF}_{AllInc}\neq\emptyset$ then ${\cal U}^{SF}_{AllContr}\neq\emptyset$, while the opposite is not true. This means that the existence of a maximally discriminant compatible scoring function in which all marginal functions contribute to the global value of an alternative does not imply the existence of a maximally discriminant compatible scoring function in which all marginal value functions are increasing. 
\end{note}

As already underlined above, the sets ${\cal U}^{SF}$, ${\cal U}^{SF}_{AllContr}$ or ${\cal U}^{SF}_{AllInc}$ could be composed of more than one maximally discriminant compatible scoring function. Each of them assigns a different value to each alternative under consideration and, therefore, a procedure to discover a well-distributed sample of maximally discriminant compatible scoring functions in these sets should be defined. In the following, we shall describe in detail how to compute a well-distributed sample of maximally discriminant compatible scoring functions in ${\cal U}^{SF}$. However, the same procedure can be easily adapted to compute well-distributed sample of maximally discriminant compatible scoring functions in ${\cal U}^{SF}_{AllContr}$ or ${\cal U}^{SF}_{AllInc}$. 

Before describing the mentioned procedure, let us observe that an additive value function $U$, as the one in (\ref{AdditiveVF}), is uniquely defined by the marginal values assigned to the evaluations $x_j^{k}$ from each marginal value function $u_{j}$, that is, $u_{j}\left(x_j^{k}\right)$ for all $g_j\in G$ and for all $k=0,\ldots,n_j$. For the sake of simplicity, denoting by $u_{j}^{k}$ the values $u_j\left(x_j^{k}\right)$, an additive value function $U$ can also be represented by the vector $\displaystyle U=\left[u_j^{k}\right]_{\substack{g_j\in G\\ k=0,\ldots,n_j}}$. 

Let as assume that ${\cal U}^{SF}\neq\emptyset$, and let us consider $U^{1}=\left[u_j^{k,1}\right]$ the maximally discriminant compatible scoring function obtained solving $LP_0$. Another maximally discriminant compatible scoring function $U\in{\cal U}^{SF}$ is different from $U^{1}$ iff $u_j^{k}\neq u_j^{k,1}$ for at least one criterion $g_j\in G$ and for at least one $k\in\{0,\ldots,n_j\}$. Therefore, the existence of a second maximally discriminant compatible value function, different from $U^1$, can be checked by solving the following MILP problem that we shall denote by MILP-1: 

$$
\begin{array}{l}
\;\;\;\;\delta_2^{*}=\displaystyle\max\delta\;\;\mbox{subject to}, \\[0,2cm]
\left.
\begin{array}{l}
\;\;\eta=\eta^{*}, \\[0,2cm]
\;\;E^{SF}, \\[0,2cm]
\;\;\delta\geqslant\delta_{min},\\[0,2cm]
\left.
\begin{array}{l}
u_j^{k}\geqslant u_{j}^{k,1}+\delta-My_{j,1}^{k,1}, \;\;k=0,\ldots,n_j,\\[0,2cm]
u_j^{k}+\delta\leqslant u_{j}^{k,1}+My_{j,2}^{k,1}, \;\;k=0,\ldots,n_j,\\[0,2cm]
y_{j,1}^{k,1},y_{j,2}^{k,1}\in\{0,1\},\\[0,2cm]
\displaystyle\displaystyle\sum_{j=1}^{m}\sum_{k=0}^{n_j}\left[y_{j,1}^{k,1}+y_{j,2}^{k,1}\right]\leqslant 2\cdot\displaystyle\sum_{g_j\in G}\left(n_j+1\right)-1
\end{array}
\right\}E_{1}
\end{array}
\right\}
\end{array}
$$ 
\noindent where
\begin{itemize}
\item $M$ is a big positive number, while $\delta$ is the minimal difference between the marginal value of the same performance $x_j^{k}$ on the maximally discriminant compatible scoring functions $U^1$ and $U$ that needs to be maximized for at least one criterion $g_j\in G$ and at least one $k=0,\ldots,n_j$. Moreover, to ensure that this difference is not too low, we fix a lower bound for this variable being $\delta_{min}$\footnote{For example, in the case study presented in Section \ref{CaseStudySection}, we fix $\delta_{min}=0.1$.}. Of course, the choice of $\delta_{min}$ will influence the computation of the other maximally discriminant compatible scoring functions obtained in addition to $U^1$ in the well-diversified sample we are looking for,   
\item $u_j^{k}\geqslant u_{j}^{k,1}+\delta-My_{j,1}^{k,1}$ translates the constraint $u_j^{k}>u_{j}^{k,1}$. In particular, if $y_{j,1}^{k,1}=1$, then the constraint is always satisfied and, therefore, $u_j^k\leqslant u_{j}^{k,1}$. If, instead, $y_{j,1}^{k,1}=0$, then the constraint is reduced to $u_j^{k}\geqslant u_{j}^{k,1}+\delta$ and, therefore, $u_j^{k}>u_{j}^{k,1}$;
\item $u_j^{k}+\delta\leqslant u_{j}^{k,1}+My_{j,2}^{k,1}$ translates the constraint $u_j^{k}<u_{j}^{k,1}$. In particular, if $y_{j,2}^{k,1}=1$, then the constraint is always satisfied and, therefore, $u_j^k\geqslant u_{j}^{k,1}$. If, instead, $y_{j,2}^{k,1}=0$, then the constraint is reduced to $u_j^{k}+\delta\leqslant u_{j}^{k,1}$ and, therefore, $u_j^{k}<u_{j}^{k,1}$;
\item the constraint $\displaystyle\sum_{j=1}^{m}\sum_{k=0}^{n_j}\left[y_{j,1}^{k,1}+y_{j,2}^{k,1}\right]\leqslant 2\cdot\displaystyle\sum_{g_j\in G}\left(n_j+1\right)-1$ is used to impose that at least one binary variable is equal to 0 and, consequently, at least one $u_j^{k}$ is such that $u_j^{k}\neq u_j^{k,1}$.
\end{itemize}

Solving $MILP-1$, two cases can occur: 
\begin{description}
\item[case 1)] $MILP-1$ is feasible: there is at least one maximally discriminant compatible scoring function in ${\cal U}^{SF}$ different from $U^1$ and the marginal values $u_j^{k}$ obtained solving $MILP-1$ define such a function. We shall denote by $U^{2}=\left[u_{j}^{k,2}\right]$ the marginal values defining the new maximally discriminant compatible scoring function;
\item[case 2)] $MILP-1$ is infeasible: the maximally discriminant compatible scoring function obtained solving $LP_{0}$ is unique and, therefore, considering $\delta_{min}$, the well-diversified sample of maximally discriminant compatible scoring functions we are looking contains only the value function $U_1$. Of course, this is true for the fixed value $\delta_{min}$ because reducing the value of $\delta_{min}$ could bring to the discovery of other maximally discriminant compatible scoring functions in addition to $U^{1}$ in the well-diversified sample we are looking for. 
\end{description}
In case 1) one can proceed in a iterative way to find a sample of well distributed maximally discriminant compatible scoring functions in ${\cal U}^{SF}$ with the considered $\delta_{min}$ in the sample. Let us assume that $t$ functions in ${\cal U}^{SF}$ have been already computed; the $(t+1)-th$ function in ${\cal U}^{SF}$ in the sample is obtained solving the following MILP problem 

$$
\begin{array}{l}
\;\;\;\delta_{t+1}^{*}=\displaystyle\max\delta\;\;\mbox{subject to}, \\[0,2cm]
\left.
\begin{array}{l}
\;\;\eta=\eta^{*}, \\[0,2cm]
\;\;E^{SF}, \\[0,2cm]
\;\;\delta\geqslant\delta_{min}, \\[0,2cm]
\;\;E_{1}\cup E_2\cup\cdots\cup E_{t}
\end{array}
\right\}
\end{array}
$$ 

\noindent and $E_{r}$ is, in general, the following set of constraints:
$$
\left.
\begin{array}{l}
u_j^{k}\geqslant u_{j}^{k,r}+\delta-My_{j,1}^{k,r}, \\[0,2cm]
u_j^{k}+\delta\leqslant u_{j}^{k,r}+My_{j,2}^{k,r}, \\[0,2cm]
y_{j,1}^{k,r},y_{j,2}^{k,r}\in\{0,1\},\\[0,2cm]
\displaystyle\displaystyle\sum_{j=1}^{m}\sum_{k=0}^{n_j}\left[y_{j,1}^{k,r}+y_{j,2}^{k,r}\right]\leqslant 2\cdot\displaystyle\sum_{g_j\in G}\left(n_j+1\right)-1.
\end{array}
\right\}E_r
$$	

Let us conclude this section observing that the sample of maximally discriminant compatible scoring functions obtained, can be used, if necessary, as a starting point for a more detailed exploration going beyond the mere sample. In fact, each convex combination $\overline{U}$ of maximally discriminant compatible scoring functions $U^1,\ldots,U^r$, that is,
$$
\overline{U}=\lambda_1 U^1 + \ldots + \lambda_r U^r, \;\;\mbox{with}\;\; \lambda_t \geqslant 0, \;\mbox{for all}\; t\in\{1,\ldots,r\},\;\;\mbox{and}\;\;\displaystyle \sum_{t=1}^{r}\lambda_t=1,
$$

\noindent is, in turn, another maximally discriminant compatible scoring function. Indeed, as $U^1, \ldots, U^r$ satisfy constraints from $E^{SF}\cup\{\eta=\eta^*\}$, also their convex combination $\overline{U}$ satisfies the same constraints. This means that starting from the obtained sample of maximally discriminant compatible scoring functions, others can be easily and meaningfully obtained, permitting the DM to examine how the ranking and the scores of alternatives change in the space of maximally discriminant compatible scoring functions.

\section{Comparison with other methods}\label{ComparisonSection}
In this section we test our scoring procedure with respect to its capacity to predict preferences on the basis of some preference comparisons provided by the DM. In this perspective, our procedure will be compared with other fourteen methods representing the state of the art in literature to obtain a ranking on the basis of the knowledge of the PWIs represented in the matrix $PWM=\left[p(a,b)\right]$. First of all, let us briefly review the other fourteen methods with which we shall compare our scoring procedure (more details on them could be found in the publications in which they have been presented). The first three methods \citep{Vetschera2017} aim to define a complete order (that is a complete, asymmetric and transitive binary relation) that optimally represents the PWIs. In particular the complete order on $A$ is represented by the 0-1 variables $y_{ab}, a,b \in A$ such that if $y_{ab}=1$, then $a$ is preferred to $b$, while this not the case if $y_{ab}=0$. The properties of completeness and asymmetry are ensured by the constraints $y_{a,b}+y_{b,a}=1, \;a,b \in A, a \neq b$, while the transitivity is obtained through the constraints $y_{a,b}\geqslant y_{a,c}+y_{c,b}-1.5, \;\mbox{for all}\; a,b,c \in A$ with $c \in A\setminus \left\{a,b\right\}$. Considering three different goodness indicators, the three following methods $\mathbf{M_1, M_2}$ and  $\mathbf{M_3}$ are then obtained.  
   
\begin{description}
\item[$\mathbf{M_{1}}:$]
$$
\begin{array}{l}
\;\max\displaystyle\sum_{\substack{(a,b)\in A\times A,\\ \;a\neq b}} p(a,b)y_{ab}, \;\mbox{subject to}\\[0,4mm]
\left.
\begin{array}{l}
y_{ab}+y_{ba}=1,\\[0,2mm]
y_{ab}\geqslant y_{ac}+y_{cb}-1.5, \;\forall c\in A\setminus\{a,b\} \\[0,2mm]
y_{ab}\in\{0,1\}
\end{array}
\right\}\forall (a,b)\in A\times A, \;a\neq b
\end{array}
$$
\item[$\mathbf{M_{2}}:$]
$$
\begin{array}{l}
\;\max\displaystyle\sum_{\substack{(a,b)\in A\times A,\\ \;a\neq b}} \log\left(p(a,b)\right)y_{ab}, \;\mbox{subject to}\\[0,4mm]
\left.
\begin{array}{l}
y_{ab}+y_{ba}=1,\\[0,2mm]
y_{ab}\geqslant y_{ac}+y_{cb}-1.5, \;\forall c\in A\setminus\{a,b\} \\[0,2mm]
y_{ab}\in\{0,1\}
\end{array}
\right\}\forall (a,b)\in A\times A, \;a\neq b
\end{array}
$$
\item[$\mathbf{M_{3}}:$]
$$
\begin{array}{l}
\;\max f_{MM}\\[0,4mm]
\left.
\begin{array}{l}
f_{MM}\leqslant p(a,b)+(1-y_{ab}),\\[0,2mm]
y_{ab}+y_{ba}=1,\\[0,2mm]
y_{ab}\geqslant y_{ac}+y_{cb}-1.5, \;\forall c\in A\setminus\{a,b\} \\[0,2mm]
y_{ab}\in\{0,1\}.
\end{array}
\right\}\forall (a,b)\in A\times A, \;a\neq b
\end{array}
$$
\end{description} 

The following methods $\mathbf{M_4-M_{15}}$ rank the alternatives from $A$ according to the increasing order of the values assigned to alternatives $a$ from $A$  by the following ranking functions: 

\begin{description}
\item[$\mathbf{M_{4}}:$] \textit{The positive outranking index}: 
$$
PosOI(a,A,PWM)=\frac{1}{|A|-1}\displaystyle\sum_{\substack{b\in A\setminus\{a\}: \\[0,1mm] p(a,b)\geqslant 0.5}}p(a,b);
$$
\item[$\mathbf{M_{5}}:$] \textit{Max in favor}: 
$$
MF(a,A,PWM)=\displaystyle\max_{b\in A\setminus\{a\}} p(a,b);
$$
\item[$\mathbf{M_{6}}:$] \textit{Min in favor}: 
$$
mF(a,A,PWM)\footnote{The function is also known as Simpson score as presented in \cite{LeskinenEtAl2006}}=\displaystyle\min_{b\in A\setminus\{a\}}p(a,b);
$$ 
\item[$\mathbf{M_{7}}:$] \textit{Sum in favor}:
$$
SF(a,A,PWM)=\displaystyle\sum_{b\in A\setminus\{a\}}p(a,b);
$$
\item[$\mathbf{M_{8}}:$] \textit{Max against}:
$$
MA(a,A,PWM)=\displaystyle-\max_{b\in A\setminus\{a\}}p(b,a);
$$
\item[$\mathbf{M_{9}}:$] \textit{Min against}:
$$
mA(a,A,PWM)=\displaystyle-\min_{b\in A\setminus\{a\}}p(b,a);
$$
\item[$\mathbf{M_{10}}:$] \textit{Sum against}:
$$
SA(a,A,PWM)=\displaystyle-\sum_{b\in A\setminus\{a\}}p(b,a);
$$
\item[$\mathbf{M_{11}}:$] \textit{Max difference}: 
$$
MD(a,A,PWM)=\displaystyle\max_{b\in A\setminus\{a\}}\left[p(a,b)-p(b,a)\right];
$$
\item[$\mathbf{M_{12}}:$] \textit{Min difference}: 
$$
mD(a,A,PWM)=\displaystyle\min_{b\in A\setminus\{a\}}\left[p(a,b)-p(b,a)\right];
$$
\item[$\mathbf{M_{13}}:$] \textit{Sum of differences}: 
$$
SD(a,A,PWM)=\displaystyle\sum_{b\in A\setminus\{a\}}\left[p(a,b)-p(b,a)\right];
$$
\item[$\mathbf{M_{14}}:$] \textit{Copeland score}: 
$$
Cop(a,A,PWM)=\displaystyle\sum_{\substack{b\in A\\ b\neq a}}r_{ab}, \;\;\;\;\mbox{where}\;\;r_{ab}=\left\{\begin{array}{ll}1, & \mbox{if}\;\;  p(a,b)\geqslant 0.5 \\[1mm] -1, & \mbox{otherwise.} \end{array}\right.
$$
\end{description}
Methods $\mathbf{M_{5}}$ - $\mathbf{M_{13}}$ have been presented in \cite{KadzinskiMichalski2016}, while methods $\mathbf{M}_4$ and $\mathbf{M}_{14}$ are reported in \cite{LeskinenEtAl2006}. Moreover, methods $\mathbf{M_4}-\mathbf{M_{14}}$ assign a unique value to each alternative in $A$ and, the greater the value assigned to the alternative, the better it is. A complete ranking of the alternatives at hand can be therefore obtained on the basis of the assigned values. In the following, we shall denote by $Ranking_{\mathbf{M}}$ the ranking obtained by the method $\mathbf{M}\in\{\mathbf{M_1},\ldots,\mathbf{M}_{14},\mathbf{ScPr}\}$, where $\mathbf{ScPr}$ will denote the scoring procedure based on the computation of the $LP_{0}$ problem we presented in Section \ref{ScoringMethodSection}.
	
\subsection{Simulations details}\label{SimulationsSection}
In order to compare $\mathbf{ScPr}$ with the fourteen methods $\mathbf{M_1-M_{14}}$ reviewed above, we shall perform an extensive set of simulations considering different decision problems in which an ``artificial" DM whose preferences are represented by a random generated value function ranks order $n$ alternatives evaluated on $m$ criteria where $n\in\{6,9,12,15\}$ and $m\in\{3,5,7\}$. Moreover, for each $(n,m)\in\{6,9,12,15\}\times\{3,5,7\}$, 500 independent runs will be done. Algorithm \ref{SimulationAlgorithm} presents the steps that have to be performed in each of the considered runs. These steps are described in the following lines:

\begin{algorithm*}
\caption{Single run steps \label{SimulationAlgorithm}}
\begin{algorithmic}
\REPEAT
\STATE 1: Generate a performance matrix of $n$ alternatives and $m$ criteria
\STATE 2: Build the DM's value function 
\STATE 3: Compute the ranking of alternatives at hand by using the artificial DM's value function and denote it by $Ranking_{\mathbf{DM}}$ 
\STATE 4: Elicit the artificial DM's preferences
\STATE 5: Sample 10,000 value functions compatible with the artificial DM's preferences and compute the PWIs 
\STATE 6: Apply method $\mathbf{M}$ with $\mathbf{M}\in\{\mathbf{M_1},\ldots,\mathbf{M}_{14},\mathbf{ScPr}\}$ to get the ranking of the alternatives at hand and denote this ranking by $Ranking_{\mathbf{M}}$ 
\STATE 7: Compute the Kendall-Tau correlation coefficient between $Ranking_{\mathbf{DM}}$ and $Ranking_{\mathbf{M}}$ for all $\mathbf{M}\in\{\mathbf{M_1},\ldots,\mathbf{M}_{14},\mathbf{ScPr}\}$
\UNTIL 500 runs have not been performed
\STATE 8: Compute the average Kendall-Tau correlation coefficient between $Ranking_{\mathbf{DM}}$ and $Ranking_{\mathbf{M}}$ for all $\mathbf{M}\in\{\mathbf{M_1},\ldots,\mathbf{M}_{14},\mathbf{ScPr}\}$
\end{algorithmic}
\end{algorithm*}

\begin{description}
\item[1:] A performance matrix $\textit{PM}=\left[pm\right]_{\substack{i=1,\ldots,n \\j=1,\ldots,m}}$ composed of $n$ rows and $m$ columns is built. The $i$-th row of the matrix $\left(pm_{i1},\ldots,pm_{im}\right)$ is a vector of $m$ values taken randomly in a uniform way in the $\left[0,1\right]$ interval and representing the evaluations of alternative $a_i\in A$ on criteria $g_j\in G$, that is, $pm_{ij}=g_j(a_i)$. Moreover, the performance matrix is built so that the alternatives having as evaluations the values in the rows of the performance matrix are non-dominated\footnote{For each $a_{i_1},a_{i_2}\in A$, $\exists g_{j_1},g_{j_2}\in G$ such that $pm_{i_1j_1}>pm_{i_2j_1}$ and $pm_{i_2j_2}<pm_{i_2j_2}$};
\item[2:] Let us assume that the artificial DM's value function is a weighted sum such that for each alternative $a_i\in A$
$$
WS(pm_{i1},\ldots,pm_{im})=WS_{i}=w_{1}pm_{i1}+\ldots+w_{m}pm_{im}=\displaystyle\sum_{j=1}^{m}w_jpm_{ij}
$$
\noindent where $w_j>0$ for all $j=1,\ldots,m$ and $\displaystyle\sum_{j=1}^{m}w_j=1$. To simulate the DM's value function, we then sample $m$ positive values $w^{DM}_1,\ldots,w^{DM}_m$ such that their sum is 1 following the procedure proposed by \cite{Rubinstein1982}\footnote{Take randomly $m-1$ values $v_1,\ldots,v_{m-1}$ in the interval $\left[0,1\right]$ and reorder the values in the set $\left\{0,v_1,\ldots,v_{m-1},1\right\}$ in a non-decreasing way so that $0=v_{(0)}\leqslant v_{(1)}\leqslant v_{(2)}\leqslant\ldots\leqslant v_{(m-1)}\leqslant v_{(m)}=1$. For each $j=1,\ldots,m$, put $w_j=v_{(j)}-v_{(j-1)}$.},
\item[3:] Apply the artificial DM's value function defined by the vector $\left(w_1^{DM},\ldots,w_m^{DM}\right)$ to compute the weighted sum of each alternative. On the basis of the values assigned to all alternatives compute their ranking and denote it by $Ranking_{\mathbf{DM}}$;
\item[4:] To simulate the artificial DM's preferences we use the procedure proposed in \cite{Vetschera2017}: we sample an alternative $a_i\in A$ and, then, we compare $a_i$ with all the other alternatives $a_{i_1}$, that is, $a_{i_1}\in A\setminus\{a_{i}\}$, so that if $WS_{i}>WS_{i_1}$, then $a_i\succ_{DM}a_{i_1}$, if $WS_{i}<WS_{i_1}$ then $a_{i_1}\succ_{DM} a_i$, while if $WS_{i}=WS_{i_1}$, then $a_i\sim_{DM}a_{i_1}$,
\item[5:] Sample 10,000 value functions compatible with the preferences provided by the artificial DM. In this case, we assume that a compatible value function is a vector of weights $(w_1,\ldots,w_m)$ so that the following set of constraints is satisfied:
$$
\left.
\begin{array}{l}
w_1pm_{i1}+\ldots+w_mpm_{im}\geqslant w_1pm_{i_11}+\ldots+w_mpm_{i_1m}+\varepsilon, \;\;\mbox{if}\;\;a_i\succ_{DM}a_{i_1}, \\[2mm]
w_1pm_{i1}+\ldots+w_mpm_{im}=w_1pm_{i_11}+\ldots+w_mpm_{i_1m}, \;\;\mbox{if}\;\;a_i\sim_{DM}a_{i_1}, \\[2mm] 
w_j\geqslant\varepsilon,\;\;\mbox{for all}\;j=1,\ldots,m,\\[2mm]
\displaystyle\sum_{j=1}^{m}w_j=1,\\[2mm]
\varepsilon>0.
\end{array}
\right\}E^{DM}_{WS}
$$
Since the artificial DM's value function is a weighted sum and the value function used to approximate its preferences is a weighted sum as well, there exists at least one compatible value function and, therefore, we can sample 10,000 compatible value functions (for example using the HAR method) from the space defined by the constraints in $E^{DM}_{WS}$. Denoting by ${\cal U}$ the set of sampled compatible value functions we have $|{\cal U}|=10,000$. Computing the sets ${\cal U}_{a\succ b}$ and ${\cal U}_{a\sim b}$ for each $(a,b)\in A\times A$ with $a\neq b$ as in eq. (\ref{SubsetsValueFunctions}), we then can compute the PWIs $p(a,b)$ as defined in eq. (\ref{PWI});
\item[6:] Apply the fourteen methods presented above as well as the scoring procedure we are proposing to compute the alternatives' ranking $Ranking_{\mathbf{M}}$ with $\mathbf{M}\in\{\mathbf{M}_1,\ldots,\mathbf{M}_{14},\mathbf{ScPr}\}$;
\item[7:] For each $\mathbf{M}\in\{\mathbf{M}_1,\ldots,\mathbf{M}_{14},\mathbf{ScPr}\}$ compute the Kendall-Tau \citep{Kendall1938} between $Ranking_{\mathbf{M}}$ and $Ranking_{\mathbf{DM}}$ and denote it as $\tau(\mathbf{M},\mathbf{DM})$;
\item[8:] For each $\mathbf{M}\in\{\mathbf{M}_1,\ldots,\mathbf{M}_{14},\mathbf{ScPr}\}$ compute the average Kendall-Tau over the 500 independent runs between $Ranking_{\mathbf{M}}$ and $Ranking_{\mathbf{DM}}$. Denote by $\overline{\tau}(\mathbf{M},\mathbf{DM})$ this value.
\end{description}
\renewcommand\arraystretch{1.3}

\begin{table}[!h]
\begin{center}
\caption{Average (over 500 independent runs) Kendall-Tau correlation coefficient between $Ranking_{\mathbf{DM}}$ and $Ranking_{\mathbf{M}}$, where $\mathbf{M}\in\{\mathbf{M}_1,\ldots,\mathbf{M}_{14},\mathbf{ScPr},\mathbf{ScPr}'\}$. In italics, the maximum average Kendall-Tau coefficient considering $\mathbf{M}\in\{\mathbf{M}_1,\ldots,\mathbf{M}_{14},\mathbf{ScPr}\}$, that is, $\overline{\tau}_{Max}=\displaystyle\max_{\mathbf{M}\in\{\mathbf{M}_1,\ldots,\mathbf{M}_{14},\mathbf{ScPr}\}}\{\overline{\tau}\left(\mathbf{M},\mathbf{DM}\right)\}$; in grey, the maximum average Kendall-Tau coefficient considering $\mathbf{M}\in\{\mathbf{M}_1,\ldots,\mathbf{M}_{14},\mathbf{ScPr},\mathbf{ScPr}'\}$, that is, $\overline{\tau}'_{Max}=\displaystyle\max_{\mathbf{M}\in\{\mathbf{M}_1,\ldots,\mathbf{M}_{14},\mathbf{ScPr},\mathbf{ScPr}'\}}\{\overline{\tau}\left(\mathbf{M},\mathbf{DM}\right)\}$. \label{TauValues}}	
\resizebox{0,9\textwidth}{!}{
		\begin{tabular}{ccccccccccccc}
		& \multicolumn{12}{c}{(\# Alternatives, \# Criteria)} \\
		\cline{2-13} 
    $\overline{\tau}({\mathbf{M}},{\mathbf{DM}})$ & \textbf(6,3) & \textbf(6,5) & \textbf(6,7) & \textbf(9,3) & \textbf(9,5) & \textbf(9,7) & \textbf(12,3) & \textbf(12,5) & \textbf(12,7) & \textbf(15,3) & \textbf(15,5) & \textbf(15,7) \\
		\hline
		$\mathbf{M}_{1}$   & 0.8288 & 0.7968 & 0.7792 & 0.8470 & 0.7828 & \textit{0.7570} & 0.8721 & 0.8107 & 0.7680 & 0.8735 & 0.8046 & 0.7832 \\
		$\mathbf{M}_{2}$   & 0.8288 & 0.7968 & 0.7792 & 0.8470 & 0.7828 & \textit{0.7570} & 0.8721 & 0.8109 & 0.7682 & 0.8735 & 0.8044 & 0.7832 \\
		$\mathbf{M}_{3}$   & 0.8288 & 0.7968 & 0.7792 & 0.8471 & 0.7828 & 0.7568 & 0.8720 & 0.8104 & 0.7682 & 0.8734 & 0.8040 & 0.7833 \\
		$\mathbf{M}_{4}$   & 0.8280 & 0.7965 & 0.7795 & \textit{0.8472} & 0.7820 & 0.7557 & \textit{0.8724} & 0.8112 & 0.7689 & 0.8737 & 0.8046 & 0.7829 \\
		$\mathbf{M}_{5}$   & 0.6185 & 0.6632 & 0.6723 & 0.5352 & 0.5916 & 0.6076 & 0.4767 & 0.5401 & 0.5737 & 0.4301 & 0.5037 & 0.5514 \\
		$\mathbf{M}_{6}$   & 0.6456 & 0.6566 & 0.6497 & 0.5678 & 0.5967 & 0.6109 & 0.5215 & 0.5614 & 0.5774 & 0.4700 & 0.5178 & 0.5391 \\
		$\mathbf{M}_{7}$   & \cellcolor{lightgray}{\textit{0.8331}} & 0.7904 & 0.7773 & 0.8466 & \textit{0.7841} & 0.7543 & 0.8703 & 0.8109 & 0.7688 & 0.8727 & 0.8041 & 0.7837 \\
		$\mathbf{M}_{8}$   & 0.6456 & 0.6566 & 0.6497 & 0.5678 & 0.5967 & 0.6109 & 0.5215 & 0.5614 & 0.5774 & 0.4699 & 0.5178 & 0.5391 \\
		$\mathbf{M}_{9}$   & 0.6185 & 0.6632 & 0.6723 & 0.5352 & 0.5916 & 0.6076 & 0.4767 & 0.5401 & 0.5737 & 0.4300 & 0.5037 & 0.5514 \\
		$\mathbf{M}_{10}$  & \cellcolor{lightgray}{\textit{0.8331}} & 0.7904 & 0.7773 & 0.8466 & \textit{0.7841} & 0.7542 & 0.8703 & 0.8109 & 0.7688 & 0.8726 & 0.8041 & 0.7837 \\
		$\mathbf{M}_{11}$  & 0.6185 & 0.6632 & 0.6723 & 0.5352 & 0.5916 & 0.6076 & 0.4767 & 0.5401 & 0.5737 & 0.4301 & 0.5037 & 0.5514 \\
		$\mathbf{M}_{12}$  & 0.6456 & 0.6566 & 0.6497 & 0.5678 & 0.5967 & 0.6109 & 0.5215 & 0.5614 & 0.5774 & 0.4700 & 0.5178 & 0.5391 \\
		$\mathbf{M}_{13}$  & \cellcolor{lightgray}{\textit{0.8331}} & 0.7904 & 0.7773 & 0.8466 & \textit{0.7841} & 0.7542 & 0.8703 & 0.8109 & 0.7688 & 0.8727 & 0.8041 & 0.7837 \\
		$\mathbf{M}_{14}$  & 0.8283 & \textit{0.7972} & \textit{0.7797} & 0.8470 & 0.7833 & 0.7569 & 0.8722 & \textit{0.8113} & \textit{0.7690} & \textit{0.8739} & \textit{0.8049} & \textit{0.7838} \\
		$\mathbf{ScPr}$    & 0.8255 & 0.7925 & 0.7742 & 0.8450 & 0.7779 & 0.7416 & 0.8699 & 0.7922 & 0.7572 & 0.8662 & 0.7893 & 0.7628 \\
		$\mathbf{ScPr}{'}$ & 0.8289 & \cellcolor{lightgray}{0.7988} & \cellcolor{lightgray}{0.7812} & \cellcolor{lightgray}{0.8476} & \cellcolor{lightgray}{0.7891} & \cellcolor{lightgray}{0.7594} & \cellcolor{lightgray}{0.8756} & \cellcolor{lightgray}{0.8186} & \cellcolor{lightgray}{0.7726} & \cellcolor{lightgray}{0.8804} & \cellcolor{lightgray}{0.8198} & \cellcolor{lightgray}{0.7948} \\
		\hline
		$deviation\left(\mathbf{ScPr}\right)$    & $-$0.91\% & $-$0.59\% & $-$0.71\% & $-$0.26\% & $-$0.79\% & $-$2.03\% & $-0.29\%$ & $-2.35\%$ & $-1.53\%$ & $-0.88\%$ & $-1.94\%$ & $-2.68\%$ \\
		$deviation\left(\mathbf{ScPr}{'}\right)$ & $-$0.50\% & $+$0.20\% & $+$0.19\% & $+$0.05\% & $+$0.64\% & $+$0.32\% & $+0.37\%$ & $+0.90\%$ & $+0.47\%$ & $+0.74\%$ & $+1.85\%$ & $+1.40\%$ \\
		\hline
		\end{tabular}
		}
\end{center}
\end{table}

For each (\# Alternatives, \# Criteria) configuration, we computed the best method between $\mathbf{M}_1$,$\ldots,$ $\mathbf{M}_{14},$ $\mathbf{ScPr}$ being the one presenting the maximum Kendall-Tau, that is, $\mathbf{\overline{M}}\in\{\mathbf{M}_1,\ldots,\mathbf{M}_{14},\mathbf{ScPr}\}$ such that
$$
\overline{\tau}\left(\mathbf{\overline{M}},\mathbf{DM}\right)=\overline{\tau}_{Max}=\displaystyle\max_{\mathbf{M}\in\{\mathbf{M}_1,\ldots,\mathbf{M}_{14},\mathbf{ScPr}\}}\{\overline{\tau}\left(\mathbf{M},\mathbf{DM}\right)\}.
$$
We put in italics in Table \ref{TauValues} the value of the average Kendall-Tau corresponding to the method presenting $\overline{\tau}_{Max}$. 

As one can see from the data in this Table, there is not a method being the best in all configurations. To have an estimate of how $\mathbf{ScPr}$ behaves with respect to the best method $\overline{\mathbf{M}}$ in each configuration, we computed the following \textit{deviation} as a ``normalized distance" of $\overline{\tau}({\mathbf{ScPr}},{\mathbf{DM}})$ from $\overline{\tau}_{Max}$ obtained as follows:

\begin{equation}\label{ErrorEquation}
deviation(\mathbf{ScPr})=\frac{\overline{\tau}({\mathbf{ScPr}},{\mathbf{DM}})-\overline{\tau}_{Max}}{\overline{\tau}_{Max}}. 
\end{equation}

\noindent The data underline that this error is in the interval $\left[-2.68\%,-0.26\%\right]$. In particular, for the configurations with three criteria, the error is always lower than $1\%$, while the configuration having the maximum error, that is, $-2.68\%$, is $(15,7)$ being the one presenting the greatest number of alternatives and criteria. However, to check if the difference between $\overline{\tau}({\mathbf{ScPr}},{\mathbf{DM}})$ and $\overline{\tau}({\overline{\mathbf{M}}},{\mathbf{DM}})$ is statistically significant we performed the Mann-Whitney U test with $5\%$ significance level on the Kendall-Tau coefficients ${\tau}({\mathbf{ScPr}},{\mathbf{DM}})$ and ${\tau}({\overline{\mathbf{M}}},{\mathbf{DM}})$ over the 500 independent runs performed for both methods.

\begin{table}[!h]
\begin{center}
\caption{Mann-Whitney U test with 5\% significance level computed on $\tau(\mathbf{\overline{M}},\mathbf{DM})$ and $\tau(\mathbf{ScPr},\mathbf{DM})$ with $\overline{\mathbf{M}}\in\{\mathbf{M}_1,\ldots,\mathbf{M}_{14},\mathbf{ScPr}\}$ such that $\overline{\tau}(\mathbf{\overline{M}},\mathbf{DM})=\displaystyle\max_{\mathbf{M}\in\{\mathbf{M}_1,\ldots,\mathbf{M}_{14},\mathbf{ScPr}\}}\{\overline{\tau}(\mathbf{M},\mathbf{DM})\}$. \label{MannWhitney}}	
\resizebox{0,9\textwidth}{!}{
		\begin{tabular}{ccccccccccccc}
		& \multicolumn{12}{c}{(\# Alternatives, \# Criteria)} \\
		\cline{2-13}
     & \textbf(6,3) & \textbf(6,5) & \textbf(6,7) & \textbf(9,3) & \textbf(9,5) & \textbf(9,7) & \textbf(12,3) & \textbf(12,5) & \textbf(12,7) & \textbf(15,3) & \textbf(15,5) & \textbf(15,7) \\
		\hline
		$p$-value & 0.9751 & 0.8573 & 0.9247 & 0.9919 & 0.9257 & 0.7186 & 0.9753 & 0.4493 & 0.6566 & 0.7959 & 0.5926 & 0.3136 \\
		\hline
		\end{tabular}
		}
\end{center}
\end{table}

As one can see from the $p$-value of the test for each (\# Alternatives, \# Criteria) configuration, the difference between the Kendall-Tau coefficients of ${\tau}({\mathbf{ScPr}},{\mathbf{DM}})$ and ${\tau}({\overline{\mathbf{M}}},{\mathbf{DM}})$ over all performed runs is not significant from the statistical point of view. Note that for all considered configurations the null hypothesis of equal medians fails to be rejected. This means that even if our scoring procedure is not the best in any configuration, it obtains an approximation of the DM's ranking of the alternatives at hand at least as good as the one produced by the best method $\overline{\mathbf{M}}$. Once more we would like to underline that, differently from the other methods with which our proposal is compared, we do not want to produce only a ranking of the alternatives under consideration. Indeed, using the results of the PWIs we want also to explain this ranking assigning to each alternative an overall evaluation (score) expressed in terms of an additive value function. Indeed, our scoring procedure builds an additive value function summarizing the PWIs that assigns an overall evaluation, that is, a score, to each alternative and that can be used to explain the contribution given by each criterion to its global value. This point will be better clarified and illustrated in the next section. 

Let us conclude this section by observing that in solving the $LP_{0}$ problem described in Section \ref{ScoringMethodSection} we can obtain $\eta^*\leqslant 0$ in some cases. Indeed, depending on the sampled compatible models on the basis of which the PWIs are computed, it is possible that there exist $a,b,c\in A$ such that $p(a,b)>0.5$ (more than 50\% of the sampled compatible models state that $a$ is at least as good as $b$), $p(b,c)>0.5$ and, however, $p(c,a)>0.5$ (similarly to what happens in the Condorcet paradox; \citealt{Condorcet1785}). Since our scoring procedure is based on the fact that the difference between $U(a)$ and $U(b)$ by the built model should be proportional to ``$p(a,b)-0.5$", we would have the following: 
\begin{description}
\item[C1)] $p(a,b)>0.5$ $\Rightarrow$ $U(a)>U(b)$,
\item[C2)] $p(b,c)>0.5$ $\Rightarrow$ $U(b)>U(c)$,
\item[C3)] $p(c,a)>0.5$ $\Rightarrow$ $U(c)>U(a)$.
\end{description}
However, constraints $\mathbf{C1)}$ and $\mathbf{C2)}$ would imply that $U(a)>U(c)$ being, of course, in contradiction with $\mathbf{C3)}$. This means that, in this case, it does not exists any compatible scoring function. Nevertheless, also in these cases, the solution of the $LP_{0}$ problem restores an additive value function that assigns a unique value to each alternative and that can be used to rank the alternatives under consideration. In fact, when $\eta^\ast \leqslant 0$, then the obtained value function $U$ is not able to perfectly represent the preferences suggested by the PWIs. However, maximizing $\eta^\ast$, that is, minimizing its absolute value since $\eta^\ast$ is non-positive, the obtained value function $U$ can be interpreted as a value function that minimally deviates from the preference order given by the PWIs. The data reported in Table \ref{TauValues} in correspondence of the $\mathbf{ScPr}$ row are obtained considering cases for which $\eta^*>0$ and cases for which $\eta^*\leqslant 0$. In Table \ref{ErrorsScPr} we reported the percentage of cases in which $\eta^*\leqslant 0$ in solving $LP_{0}$ that, as can be observed, increases with the number of alternatives or criteria considered in the performed simulations.  

\begin{table}[!h]
\begin{center}
\caption{Percentage of runs for which the solution of the $LP_{0}$ problem gave $\eta^*\leqslant 0$. \label{ErrorsScPr}}	
\resizebox{0,9\textwidth}{!}{
		\begin{tabular}{cccccccccccc}
		\multicolumn{12}{c}{(\# Alternatives, \# Criteria)} \\
		\hline
     \textbf(6,3) & \textbf(6,5) & \textbf(6,7) & \textbf(9,3) & \textbf(9,5) & \textbf(9,7) & \textbf(12,3) & \textbf(12,5) & \textbf(12,7) & \textbf(15,3) & \textbf(15,5) & \textbf(15,7) \\
		\hline
		0.4\% & 1\% & 0.8\% & 0.6\% & 2.8\% & 2.6\% & 1.4\% & 5.6\% & 3.6\% & 3.2\% & 8.2\% & 9.2\% \\
		\hline
		\end{tabular}
		}
\end{center}
\end{table}

Even if, as underlined above, the procedure we are proposing is able to build an additive value function summarizing the PWIs, we reported in Table \ref{TauValues} the average Kendall-Tau correlation coefficients computed for the method $\mathbf{ScPr}$ only for the runs in which solving $LP_{0}$ we got $\eta^*>0$. These values are reported in the row corresponding to $\mathbf{ScPr}{'}$. In the table we put in grey the maximum average Kendall-Tau correlation coefficient obtained considering now also the method $\mathbf{ScPr}{'}$, that is, the method $\overline{\overline{\mathbf{M}}}\in\{\mathbf{M}_{1},\ldots,\mathbf{M}_{14},\mathbf{ScPr},\mathbf{ScPr}'\}$ such that 

$$
\overline{\tau}(\overline{\overline{\mathbf{M}}},\mathbf{DM})=\overline{\tau}'=\max_{\mathbf{M}\in\{\mathbf{M}_1,\ldots,\mathbf{M}_{14},\mathbf{ScPr},\mathbf{ScPr}'\}}\{\overline{\tau}\left(\mathbf{M},\mathbf{DM}\right)\}.
$$

As one can see from the results shown in Table \ref{TauValues} in correspondence of the $\mathbf{ScPr}'$ method, our scoring procedure, cleaned of the runs in which $\eta^*\leqslant 0$, presents the maximum average Kendall-Tau correlation coefficient in eleven of the twelve considered cases. Moreover, in the unique case in which $\mathbf{ScPr}'$ is not the best among the considered methods, that is $(6,3)$, the error computed for $\mathbf{ScPr}'$ using the equation (\ref{ErrorEquation}) is very small, more precisely $-$0.5\%. Summarizing, we can state that: (i) considering only the cases in which the solution of $LP_{0}$ gives $\eta^*>0$, our scoring procedure is the best among the considered methods in almost all (\# Alternatives, \# Criteria) configurations and, (ii) considering all 500 runs, the results obtained by our scoring procedure in replying the artificial DM's ranking are comparable with the results obtained by the best method among the considered ones. 

\section{Case study}\label{CaseStudySection}
In this section, we shall apply the method described in the Section \ref{ScoringSection} showing its main characteristics and comparing its performances with another methodology used in literature to deal with the same type of problem. For this reason, we shall consider a financial problem in which the returns of seven funds are evaluated with respect to five performance measures (for a standard survey on performance measures see, for example, \cite{amenc2005portfolio,bacon2012practical}). The considered funds belong to the sub-class of the balanced bond funds which may invest in stocks a proportion of their assets between 10\% and 50\%. The historical data used to estimate the relevant statistics are the daily logarithmic return of the interval between the 01/01/2018 - 01/01/2021 (784 data points). The performance measure considered are: 
\begin{itemize}
	\item ${g}_{1}$: the Sharpe Ratio (SR; \citealt{sharpe1998sharpe}), defined as the ratio between the expected return and the standard deviation. It is often chosen by practitioners to rank managed portfolios and belongs to the class of reward-to-variability measures;
	\item ${g}_{2}$: the Treynor Ratio (TR; \citealt{scholz2005investor}), defined as the ratio between the expected return and the systematic risk sensitivity with respect to a benchmark (here we consider the MSCI World index \citep{bacmann2003alternative}). It assumes the form of a reward-to-variability measures but it derives from the Capital Asset Pricing Model (CAPM) portfolio theory;
	\item ${g}_{3}$: the Average Vaue-at-Risk Ratio (AVaRR; \citealt{rockafellar2000optimization}), defined as a relative performance measure of the same class of SR and TR. It looks at the average amount of potential losses suffered by the portfolio manager;
	\item ${g}_{4}$: the Jensen Alpha (JA; \citealt{jensen1968performance}), derived from the CAPM theory, defined as the excess return of a fund over the theoretical benchmark (MSCI World index); 	
	\item ${g}_{5}$: the Morningstar Risk-Adjusted Return (MRAR; \citealt{sharpe1998morningstar}), derived within the Expected Utility theory, defined as the annualized geometric average of the returns. 
\end{itemize}

The evaluations of the alternatives on the considered criteria are given in Table \ref{TablePerformances}.

\begin{table}[!h]
\begin{center}
\caption{Evaluation of the 7 funds on the five considered criteria\label{TablePerformances}}	
		\begin{tabular}{llccccc}
     & & ${g}_{1}(\cdot)$ & ${g}_{2}(\cdot)$ & ${g}_{3}(\cdot)$ & ${g}_{4}(\cdot)$ & ${g}_{5}(\cdot)$ \\
		\hline
    $a_{1}$ & Allianz Multipartner Multi20 &	0.0403 & 0,0010	& -0.0155	& -0.0030	& -0,0010  \\
    $a_{2}$ & Amundi Bilanciato Euro C	   &  0.0257 & 0.0004	& -0.0103	& -0.0014	& -0.0008  \\
    $a_{3}$ & Arca Te - Titoli Esteri	     &  0.0322 & 0.0009 &	-0.0133 & -0.0022 &	-0.0011  \\
    $a_{4}$ & Bancoposta Mix 2 A Cap       &	0.0193 & 0.0003 &	-0.0080 &	-0.0011 &	-0.0006  \\
    $a_{5}$ & Etica Rendita Bilanciata I	 &  0.0334 & 0.0005 &	-0.0150 &	-0.0009 &	-0.0013  \\
    $a_{6}$ & Eurizon Pir Italia 30 I	     &  0.0219 & 0.0003 &	-0.0088 & -0.0011 &	-0.0007  \\
    $a_{7}$ & Pramerica Global Multiasset 30 & -0.0018 & 0.0000 &	0.0007 & -0.0010 & 0.0006 \\	
		\hline
		\end{tabular}
\end{center}
\end{table}

The evaluations of the considered alternatives on the criteria at hand will be aggregated by the weighted sum (\ref{WeightedSum})
\begin{equation}\label{WeightedSum}
WS(a)=WS\left(g_1(a),\ldots,g_m(a)\right)=\displaystyle\sum_{j=1}^{m}w_j\cdot {g}_{j}(a)
\end{equation}
where $w_j$ are the weights of criteria $g_j$ and they are such that $w_j>0$ for all $g_j\in G$ and $\displaystyle\sum_{j=1}^{m}w_j=1$.

However, the use of a weighted sum implies that the evaluations are expressed on the same scale. For such a reason, before applying the weighted sum, we used a standardization technique applied in \cite{GrecoEtAl2018} transforming the evaluation of each $a\in A$ on $g_j\in G$, that is $g_j(a)$, into the standardized value $\overline{g}_j(a)$ as follows: 

$$
\overline{g}_j(a)=\left\{
\begin{array}{lll}
0 & \mbox{if} & g_j(a)\leqslant M_j-3s_j,\\[0,1mm]
0.5+\frac{g_j^{z}(a)}{6} & \mbox{if} & M_j-3s_j<g_j(a)<M_j+3s_j,\\[0,1mm]
1 & \mbox{if} & g_j(a)\geqslant M_j+s_j
\end{array}
\right.
$$
\noindent if $g_j$ has an increasing direction of preference or 
$$
\overline{g}_j(a)=\left\{
\begin{array}{lll}
0 & \mbox{if} & g_j(a)\geqslant M_j+3s_j,\\[0,1mm]
0.5-\frac{g_j^{z}(a)}{6} & \mbox{if} & M_j-3s_j<g_j(a)<M_j+3s_j,\\[0,1mm]
1 & \mbox{if} & g_j(a)\leqslant M_j-s_j
\end{array}
\right.
$$
\noindent if $g_j$ has a decreasing direction of preference. In both cases, $M_j=\frac{1}{|A|}\displaystyle\sum_{a\in A}g_j(a)$, $s_j=\sqrt{\frac{\displaystyle\sum_{a\in A}(g_j(a)-M_j)^2}{|A|}}$, and $g_j^{z}(a)=\frac{g_j(a)-M_{j}}{s_j}$. In this way, observing that the five considered criteria have an increasing direction of preference, the values in Table \ref{TablePerformances} are transformed in those given in Table \ref{TableNormalizedPerformances}. 

\begin{table}[!h]
\begin{center}
\caption{Normalized values of the 7 funds on the five criteria at hand\label{TableNormalizedPerformances}}	
		\begin{tabular}{llccccc}
     & & $\overline{g}_{1}(\cdot)$ & $\overline{g}_{2}(\cdot)$ & $\overline{g}_{3}(\cdot)$ & $\overline{g}_{4}(\cdot)$ & $\overline{g}_{5}(\cdot)$ \\
		\hline
    $a_{1}$ & Allianz Multipartner Multi20 &	0.6940 &	0.7349	& 0.3370 &	0.1917 &	0.4171 \\
    $a_{2}$ & Amundi Bilanciato Euro C	   &  0.5157 &	0.4869  &	0.4924 &	0.5268 &	0.4651 \\
    $a_{3}$ & Arca Te - Titoli Esteri	     &  0.5943 &	0.6939	& 0.4013 &	0.3492 &	0.3923 \\
    $a_{4}$ & Bancoposta Mix 2 A Cap       &	0.4370 &  0.3980  &	0.5596 &	0.5880 &	0.5381 \\
    $a_{5}$ & Etica Rendita Bilanciata I	 &  0.6102 &	0.4950  &	0.3522 &	0.6333 &	0.3455 \\
    $a_{6}$ & Eurizon Pir Italia 30 I	     &  0.4694 &	0.4342  &	0.5357 &	0.5926 &	0.4933 \\
    $a_{7}$ & Pramerica Global Multiasset 30	& 0.1793 &	0.2571 &	0.8219 &	0.6184 &	0.8487 \\
		\hline
		\end{tabular}
\end{center}
\end{table}

For the sake of simplicity let us assume that there is not any preference information provided by the DM and, consequently, the space from which the compatible models (the weighted sum in our case) have to be sampled is the following: 
$$
W=\{(w_1,\ldots,w_5)\in\rea^{5}: w_j>0, \;\forall j=1,\ldots,5,\;\mbox{and}\;\displaystyle\sum_{j=1}^{5}w_j=1\}.
$$

Sampling 100,000 weight vectors $\left(w_1,\ldots,w_5\right)$ from the space $W$ and computing the value assigned by the weighted sum to the seven funds for each of the 100,000 weight vectors, we obtain the PWIs shown in Table \ref{PWITable} being the basis of our scoring procedure that builds an additive value function (see eq. \ref{AdditiveVF}) summarizing the PWIs given in Table \ref{PWITable}.

\begin{table}[!h]
\begin{center}
\caption{Pairwise winning indices of the considered funds \label{PWITable}}	
		\begin{tabular}{cccccccc}
    $p(\cdot,\cdot)$ & $a_{1}$ & $a_{2}$ & $a_{3}$ & $a_{4}$ & $a_{5}$ & $a_{6}$ & $a_{7}$ \\
		\hline
    $a_{1}$ & 0.50 &	0.4166	& 0.4261	& 0.4082 &	0.5237 &	0.4070 &	0.3673 \\
    $a_{2}$ & 0.5834 &	0.50	& 0.5928 &	0.4074 &	0.5577 &	0.3673 &	0.3681 \\
    $a_{3}$ & 0.5739 &	0.4072 &	0.50	& 0.4050 &	0.5482 &	0.3969 &	0.3729 \\
    $a_{4}$ & 0.5918 &	0.5926 &	0.5950 &	0.50	& 0.5595 &	0.4570 &	0.3586 \\
    $a_{5}$ & 0.4763 &	0.4423 &	0.4519 &	0.4405 &	0.50	& 0.4145 &	0.3933 \\
    $a_{6}$ & 0.5930 &	0.6327 &	0.6031 &	0.5431 &	0.5855 &	0.50	& 0.3816 \\
    $a_{7}$ & 0.6327 &	0.6319 &	0.6271 &	0.6414 &	0.6067 &	0.6184 &	0.50 \\
		\hline
		\end{tabular}
\end{center}
\end{table}

Following the notation introduced in Section \ref{ScoringSection}, since on each criterion all alternatives have different evaluations, $n_{j}=6$, for all $j=1,\ldots,5$. Consequently, an additive value function is defined by five marginal value functions assigning values to seven different evaluations, that is, $U=\left[u_j^k\right]_{\substack{j=1,\ldots,5 \\ k=0,\ldots,6}}.$

\begin{figure}
\centering
\subfigure[AVaRR ($\overline{g}_3$)\label{g3}]
{\includegraphics[scale=0.98]{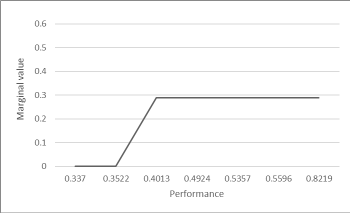}}
\subfigure[JA ($\overline{g}_4$)\label{g4}]
{\includegraphics[scale=0.98]{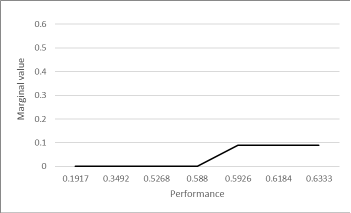}}
\subfigure[MRAR ($\overline{g}_5$)\label{g5}]
{\includegraphics[scale=0.98]{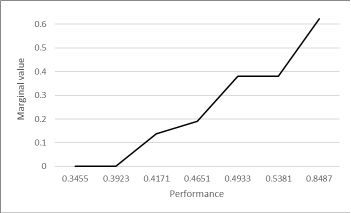}}
\caption{Marginal value functions for criteria AVaRR, JA and MRAR obtained solving $LP_0$\label{MarginalValueFunctionsLP0}}
\end{figure}

Solving the $LP_{0}$ problem presented in Section \ref{ScoringMethodSection}, we find that $E^{SF}$ is feasible and $\eta^{*}=2.0513$. Therefore, at least one compatible scoring function exists and the one obtained solving $LP_{0}$, denoted by $U^{1}$, is given in Table \ref{LP0Function}, while the respective marginal value functions are shown in Figure \ref{MarginalValueFunctionsLP0}. 

\begin{table}[!h]
\begin{center}
\caption{The additive value function obtained solving the $LP_{0}$ problem\label{LP0Function}}	
\resizebox{1\textwidth}{!}{
		\begin{tabular}{ll|ll|ll|ll|ll}
     \multicolumn{2}{c}{SR ($\overline{g}_{1}$)} & \multicolumn{2}{c}{TR ($\overline{g}_{2}$)} & \multicolumn{2}{c}{AVaRR ($\overline{g}_{3}$)} & \multicolumn{2}{c}{JA ($\overline{g}_{4}$)} & \multicolumn{2}{c}{MRAR ($\overline{g}_{5}$)} \\
		\hline
    $x_{1}^{0}=0.1793$ & $u_{1}^{0,1}=0$ & $x_{2}^{0}=0.2571$ & $u_{2}^{0,1}=0$ & $x_{3}^{0}=0.337$ & $u_{3}^{0,1}=0$ & $x_{4}^{0}=0.1917$ & $u_{4}^{0,1}=0$ & $x_{5}^{0}=0.3455$ & $u_{5}^{0,1}=0$ \\ 
		$x_{1}^{1}=0.437$ & $u_{1}^{1,1}=0$ & $x_{2}^{1}=0.398$ & $u_{2}^{1,1}=0$ & $x_{3}^{1}=0.3522$ & $u_{3}^{1,1}=0$ & $x_{4}^{1}=0.3492$ & $u_{4}^{1,1}=0$ & $x_{5}^{1}=0.3923$ & $u_{5}^{1,1}=0$ \\ 
		$x_{1}^{2}=0.4694$ & $u_{1}^{2,1}=0$ & $x_{2}^{2}=0.4342$ & $u_{2}^{2,1}=0$ & $x_{3}^{2}=0.4013$ & $u_{3}^{2,1}=0.2885$ & $x_{4}^{2}=0.5268$ & $u_{4}^{2,1}=0$ & $x_{5}^{2}=0.4171$ & $u_{5}^{2,1}=0.137$ \\ 
  	$x_{1}^{3}=0.5157$ & $u_{1}^{3,1}=0$ & $x_{2}^{3}=0.4869$ & $u_{2}^{3,1}=0$ & $x_{3}^{3}=0.4924$ & $u_{3}^{3,1}=0.2885$ & $x_{4}^{3}=0.588$ & $u_{4}^{3,1}=0$ & $x_{5}^{3}=0.4651$ & $u_{5}^{3,1}=0.1904$ \\ 
		$x_{1}^{4}=0.5943$ & $u_{1}^{4,1}=0$ & $x_{2}^{4}=0.495$ & $u_{2}^{4,1}=0$ & $x_{3}^{4}=0.5357$ & $u_{3}^{4,1}=0.2885$ & $x_{4}^{4}=0.5926$ & $u_{4}^{4,1}=0.0883$ & $x_{5}^{4}=0.4933$ & $u_{5}^{4,1}=0.3803$ \\ 
		$x_{1}^{5}=0.6102$ & $u_{1}^{5,1}=0$ & $x_{2}^{5}=0.6939$ & $u_{2}^{5,1}=0$ & $x_{3}^{5}=0.5596$ & $u_{3}^{5,1}=0.2885$ & $x_{4}^{5}=0.6184$ & $u_{4}^{5,1}=0.0883$ & $x_{5}^{5}=0.5381$ & $u_{5}^{5,1}=0.3803$ \\ 
		$x_{1}^{6}=0.694$ & $u_{1}^{6,1}=0$ & $x_{2}^{6}=0.7349$ & $u_{2}^{6,1}=0$ & $x_{3}^{6}=0.8219$ & $u_{3}^{6,1}=0.2885$ & $x_{4}^{6}=0.6333$ & $u_{4}^{6,1}=0.0883$ & $x_{5}^{6}=0.8487$ & $u_{5}^{6,1}=0.6232$ \\ 
		\hline
		\end{tabular}
		}
\end{center}
\end{table}

\noindent As one can see from the table, the greatest values on criteria SR and TR have a null marginal value meaning that the obtained maximally discriminant compatible scoring function does not assign any ``importance" to these criteria that, consequently, do not contribute to the global value obtained by each alternative. Going to the other three criteria, on the one hand, the greatest maximal shares correspond to MRAR (0.6232), followed by AVaRR (0.2885), while, on the other hand, the least maximal share is due to JA (0.0833). Looking at the changes in the marginal values, they are different for the three criteria. In particular, the greatest variation is observed for MRAR in passing from 0.0.5381 ($u_5(0.5381)=0.3803$) to 0.8487 ($u_5(0.8487)=0.6232$). Going to AVaRR and JA, the corresponding marginal value functions present only a single change in passing, on the one hand, from 0.3522 ($u_3(0.3522)=0$) to 0.4013 ($u_3(0.4013)=0.2885$) and, on the other hand, in passing from 0.588 ($u_4(0.588)=0$) to 0.5926 ($u_4(0.5926)=0.0883$). This means that changes in the performances on these two criteria affect very marginally the global value of the considered funds. 

\begin{table}[!h]
\begin{center}
\caption{Global utility of the seven funds applying the maximally discriminant compatible scoring function $U^{1}$ obtained solving $LP_{0}$\label{GlobalValueTable}}	
\resizebox{0.7\textwidth}{!}{
		\begin{tabular}{cccccccc}
     & $a_{1}$ & $a_{2}$ & $a_{3}$ & $a_{4}$ & $a_{5}$ & $a_{6}$ & $a_{7}$ \\
		\hline
		$U^{1}(\cdot)$ & 0.137 & 0.4789 & 0.2885 & 0.6688 & 0.0883 & 0.7571 & 1      \\
		\end{tabular}
		}
\end{center}
\end{table}

\begin{figure}
\centering
\includegraphics[scale=1.5]{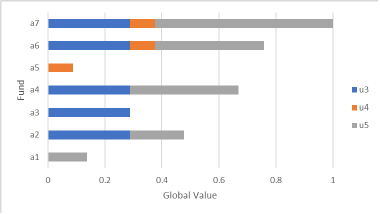}
\caption{Global value assigned to each fund by the scoring function obtained solving $LP_0$ and shown in Table \ref{LP0Function}\label{GlobalValueLP0}}
\end{figure}

The maximally discriminant compatible scoring function obtained solving $LP_{0}$ and shown in Table \ref{LP0Function} assigns a unique value to each fund as shown in Table \ref{GlobalValueTable}. The ranking of the considered funds on the basis of this function is the following:

$$
a_{7} \succ a_{6} \succ a_{4} \succ a_{2} \succ a_{3} \succ a_{1} \succ a_{5}.
$$

To better explain why each fund fills a certain position in the provided ranking and its own score, that is, its global value, is assigned we examine the contribution given to the global value of the same fund by the criteria under consideration as shown in Figure \ref{GlobalValueLP0}. For example, even if AVaRR and JA contribute in the same way to the global value of $a_6$ and $a_7$, the better global value obtained by $a_7$ can be explained by its better performance on MRAR ($\overline{g}_{5}(a_7)=0.8487$ and $\overline{g}_{5}(a_6)=0.4933$) to which corresponds a difference in the marginal value of 0.2429 ($u_{5}(0.8487)-u_{5}(0.4933)=0.6232-0.3803=0.2429$). Analogously, even if $a_5$ is better than $a_1$ on AVaRR and JA, this is not enough to compensate the greater value obtained by $a_1$ on MRAR. Indeed, for both funds only one criterion gives a marginal contribution to their global value (MRAR for $a_1$ and JA for $a_5$) since all the other performances on the remaining criteria are associated with a null marginal value. However, the marginal value given by MRAR to $a_1$ (0.137) is greater than the one given by JA to $a_5$ (0.0833) that presents the greatest performance on this criterion (0.6333). Analogous considerations can be done for the funds put in the middle of the obtained ranking, that is, $a_2$, $a_3$ and $a_4$. 

\subsection{Maximally discriminant compatible scoring functions with different characteristics and checking for other maximally discriminant compatible scoring functions}\label{ImprovementsSection}
As already described above, solving the $LP_{0}$ problem we found that ${\cal U}^{SF}\neq\emptyset$ since $\eta^{*}>0$ and, therefore, there exists at least one compatible scoring function able to summarize the information of the PWIs shown in Table \ref{PWITable}. Moreover, we observed that for the function $U^{1}$ obtained solving the LP problem and shown in Table \ref{LP0Function}, only the last three criteria (AVaRR, JA and MRAR) are giving a contribution to the global value of each fund while the marginal value attached to the best performances on SR and TR, that is $0.694$ and $0.7349$, respectively, is zero. For this reason, following the procedure shown in Section \ref{ScoringSection}, one can wonder if there is a maximally discriminant compatible scoring function such that all criteria give a contribution to the global value of the funds at hand. Solving $LP_{1}$, we find that $E_{AllContr}^{SF}$ is feasible and $h^{*}=0$. This means that each function in ${\cal U}^{SF}$ is such that at least one marginal value function gives a null contribution to the global value assigned to the alternatives and, consequently, ${\cal U}^{SF}_{AllContr}=\emptyset$. Because of the Note \ref{InclusionNote}, we find that also ${\cal U}^{SF}_{AllInc}=\emptyset$ and, therefore, there is not any maximally discriminant compatible scoring function such that all marginal value functions are monotone in their domain. 

\begin{figure}
\centering
\subfigure[AVaRR ($\overline{g}_3$)\label{Furthg3}]
{\includegraphics[scale=1.40]{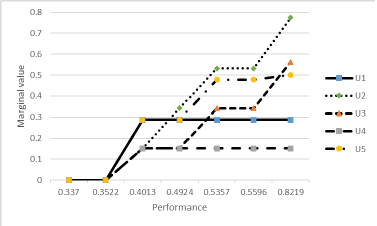}}
\subfigure[MRAR ($\overline{g}_5$)\label{Furthg5}]
{\includegraphics[scale=1.40]{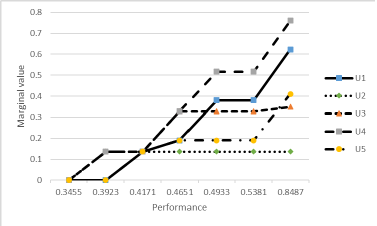}}
\caption{Marginal value functions for criteria AVaRR and MRAR for five maximally discriminant compatible value functions being the most distant among them\label{FurthestMarginalValueFunctions}}
\end{figure}

Since ${\cal U}^{SF}\neq\emptyset$, one can wonder if there exists another maximally discriminant compatible scoring function different from $U^{1}$. For this reason, we iteratively solved the MILP problems described in Section \ref{ScoringSection} considering $\delta_{min}=0.1$. In this way, in addition to the function found solving $LP_{0}$, we get a sample of twenty well-diversified maximally discriminant compatible scoring functions. An interesting aspect is that for all these functions, the marginal value corresponding to the greatest performance on criteria SR and TR is zero meaning that both criteria have a null impact on the global value of the alternatives. Moreover, the marginal value function for criterion JA is the same for the twenty functions and, in particular, is the one shown in Figure \ref{g4}. Consequently, the twenty considered maximally discriminant compatible scoring functions in the sample, together with the one obtained solving $LP_{0}$, differ only for the marginal value functions of AVaRR and MRAR. For this reason, in Figure \ref{FurthestMarginalValueFunctions} we have shown the marginal value functions with respect to these two criteria for $LP_{0}$ and other four ``most distant functions" among the other twenty obtained by the procedure described in Section \ref{ExtensionsSection}. The four functions are chosen in the following way: Let us denote by ${\cal U}_{\cal F}$ the set composed of the twenty maximally discriminant compatible scoring functions in the sample and by $U^1$ the function obtained solving $LP_{0}$. At first, let us select the maximally discriminant compatible scoring function in ${\cal U}_{\cal F}$ being the farthest from $U^1$, that is $U^2\in{\cal U}_{\cal F}$ such that $d(U^2,U^1)=\displaystyle\max_{U\in{\cal U}_{\cal F}}d(U,U^1)$, where $d(U,U^1)$ is the euclidean vector-to-vector distance between $U$ and $U^1$. Let us denote by ${\cal CF}$ the set composed of the chosen functions up to now, that is, ${\cal CF}=\{U^1,U^2\}$. After that, let us add to ${\cal CF}$ the function $U^k\in{\cal U}_{\cal F}$ presenting the maximum $\displaystyle\min_{U\in {\cal CF}}\{d(U^k,U)\}$, that is, the function presenting the maximal minimal distance from the functions that have already been included in ${\cal CF}$. The procedure continues then iteratively until $|{\cal CF}|=5$.  
	
Looking at the marginal value functions corresponding to AVaRR (Figure \ref{Furthg3}) and MRAR (Figure \ref{Furthg5}), one can observe that the five considered maximally discriminant compatible scoring functions are quite different with respect to the importance assigned to two mentioned criteria, that is, to the marginal value assigned to the greatest performance on the considered criterion. For example, in $U^{2}$ and $U^{4}$, the two marginal value functions contribute in a completely different way to the global value of the seven alternatives. On the one hand, AVaRR slightly contributes to $U^{4}$ since, for this function, $u_{3}^{4}(0.8219)=0.1515$, while it has a very great importance in $U^{2}$ since $u_{3}^{2}(0.8219)=0.7747$. On the other hand, the opposite behavior can be observed for $U^{2}$ and $U^{4}$ for criterion MRAR. Indeed, this criterion slightly contributes to the global value of the alternatives in $U^2$ since $u_{5}^{2}(0.8487)=0.1369$, while it contributes in a considerable way in $U^4$ since $u_{5}^{4}(0.8487)=0.7601$. This sheds light on the importance of taking into account not only one maximally discriminant compatible scoring function but all the maximally discriminant compatible scoring functions in the well distributed sample obtained through the iterative procedure described in Section \ref{ExtensionsSection}. 

\section{Conclusions}\label{ConclusionsSection}
In this paper we have proposed a new scoring procedure. The procedure assigns a value to each alternative under consideration summarizing the Pairwise Winning Indices (PWIs) provided by the Stochastic Multicriteria Acceptability Analysis (SMAA;\citealt{LahdelmaHokkanenSalminen1998,PelissariEtAl2020}). The method builds an additive value function compatible with the PWIs that assigns a score to each alternative and that, for this reason, is called \textit{compatible scoring function}. The idea under the proposal is that the difference between the scores assigned to two alternatives $a$ and $b$ from the built value function should be proportional to $p(a,b)$, that is, to the probability with which $a$ is considered at least as good as $b$. The greater $p(a,b)$, the larger the difference of the scores attributed to $a$ and $b$. The compatible scoring function is obtained solving a simple LP problem. In case more than one compatible scoring function maximally discriminating among the alternatives exists, we show an iterative procedure aiming to find a well-diversified sample of them. Moreover, some LP problems are presented to discover, among the maximally discriminant compatible scoring functions, some of them presenting particular characteristics: (i) functions for which the marginal value functions are strictly monotone, or (ii) functions for which all criteria contribute to the global score assigned to the alternatives from the built function. 

Several other methods summarizing the information contained in the PWIs have been proposed in literature before. However, differently from our proposal, they aim only to rank order all of them from the best to the worst. The ranking procedures proposed in the literature very often assign a value to each alternative that has not any particular meaning. Differently from all these methods our proposal has the advantage that the score assigned to each alternative is based on the construction of an additive value function that, through its marginal value functions, allows the DM to get some explanations on the reasons for which  an alternative fills a given position and obtains its specific score. In particular, the ranking position and the score of each alternative can be explained in terms of the contributions given by each criterion to the global value assigned to the considered alternative. This is a great advance from the decision aiding point of view since the scoring function we are proposing answers to the explainability concerns being nowadays very relevant for any decision aiding method (see, e.g. \citealt{ArrietaEtAl2020}). 

To prove that the new proposal, beyond explaining the rank position and the score of the alternatives, is efficient in predicting the preferences of the DM on the basis of the PWIs based on the DM's partial preference information, we performed a large set of computational experiments. We compared our scoring procedure to other fourteen methods that have been proposed inliterature and that represent the state of the art in this field. We simulated an artificial DM in problems composed of different numbers of alternatives and criteria trying to replicate the ranking produced by the artificial DM itself. We considered 6, 9, 12 and 15 alternatives and 3, 5 and 7 criteria. For each (\# Alternatives, \# Criteria) configuration we performed 500 independent runs applying the fourteen mentioned methods and our scoring procedure. To check how efficient the methods are in replicating the ranking of the artificial DM, we computed the Kendall-Tau between the preference ranking of the artificial DM and the ranking produced by each considered method. These Kendall-Tau values are then averaged over the 500 independent runs. \\
The results show that even if our proposal is not getting the best value (the maximum average Kendall-Tau) for any of the considered configurations, the ``deviation" from the best average Kendall-Tau value is always lower than 3\%. To check if the difference from the Kendall-Tau values of our proposal and the ones of the best method for each (\# Alternatives, \# Criteria) configuration is significant from the statistical point of view, we performed a Mann-Whitney U test with 5\% significance level. The test shows that in none of the considered configurations the difference is significant. This means that our scoring procedure is able to reproduce the preferences of the artificial DM and, at the same time, differently from all the other methods, it is able to give an explaination of the reasons giving to the alternatives a certain rank position. 

Finally, we have shown how to apply the new scoring procedure to a financial problem in which seven funds are evaluated with respect to five different criteria underlying the potentialities of the proposed method in explaining how the criteria contribute to the global value assigned to the alternatives. 

As further directions of research we plan to apply the new proposal to some real world decision problems to which SMAA has been applied and for which a final ranking of the alternatives under consideration has to be produced. Moreover, how to extend the scoring procedure to summarize the PWIs obtained in a problem presenting a hierarchical structure of criteria (see, for example, \citealt{CorrenteEtAl2017b}) deserves to be investigated.

\section*{Acknowledgments}
\noindent The authors wish to acknowledge the support of the Ministero dell'Istruzione, dell'Universit\'{a} e della Ricerca (MIUR) - PRIN 2017, project "Multiple Criteria Decision Analysis and Multiple Criteria Decision Theory", grant 2017CY2NCA. Moreover, Salvatore Corrente wishes to acknowledge the support of the STARTING GRANT of the University of Catania.

\bibliographystyle{plainnat}
\bibliography{Full_bibliography}

\end{document}